\documentclass{article}

\usepackage{epsfig}
\usepackage{amsmath, amsthm, amssymb}
\usepackage[all]{xy}
\input xy
\xyoption{all}
\xyoption{2cell}
\xyoption{v2}

\usepackage{hyperref}

\newcommand{\hepth}[1]{\href{http://arxiv.org/abs/hep-th/#1}{{\tt arXiv:hep-th/#1}}}

\newcommand{\Math}[1]{\href{http://arxiv.org/abs/math/#1}{{\tt arXiv:math/#1}}}

\newcommand{\arxiv}[1]{\href{http://arxiv.org/abs/#1}{{\tt arXiv:#1}}}

\newcommand{\maps}{\colon}
\newcommand{\C}{\mathcal{C}} 
\newcommand{\G}{\mathcal{G}}

\newcommand{\Aut}{\mathit{Aut}} 
\newcommand{\cU}{\mathcal{U}} 
\newcommand{\cV}{\mathcal{V}} 

\newcommand{\Ob}{\mathrm{Ob}} 
\newcommand{\Mor}{\mathrm{Mor}}

\newcommand{\cG}{\mathcal{G}} 
\newcommand{\AUT}{\mathrm{AUT}} 

\newcommand{\ZZ}{\mathbb{Z}} 
\newcommand{\QQ}{\mathbb{Q}}

\newcommand{\Spin}{\mathrm{Spin}} 
\newcommand{\String}{\mathrm{String}}
\newcommand{\U}{\mathrm{U}}

\renewcommand{\a}{\alpha} 
\renewcommand{\b}{\beta} 
\renewcommand{\c}{\gamma}

\newtheorem{theorem}{Theorem}

\newtheorem{lemma}{Lemma} 
\newtheorem{corollary}{Corollary} 

\theoremstyle{definition}
\newtheorem{definition}[theorem]{Definition}

\title{The Classifying Space of a Topological 2-Group} 
\author{John C.\ Baez\thanks{Department of Mathematics, 
University of California, 
Riverside, CA 92521, USA.
Email: \texttt{baez@math.ucr.edu}}
\and
Danny Stevenson\thanks{Fachbereich Mathematik,
Universit{\"a}t Hamburg, Hamburg, 20146, Germany,  
Email: \texttt{stevenson@math.uni-hamburg.de}} 
}

\begin{document} 
\maketitle 

\begin{abstract}
\noindent
Categorifying the concept of topological group, one obtains
the notion of a `topological 2-group'.  This in turn allows
a theory of `principal 2-bundles' generalizing the usual theory 
of principal bundles.  It is well-known that under mild conditions 
on a topological group $G$ and a space $M$, principal $G$-bundles 
over $M$ are classified by either the \v{C}ech cohomology 
$\check{H}^1(M,G)$ or the set of homotopy classes $[M,BG]$, where $BG$ 
is the classifying space of $G$.  Here we review work by Bartels, 
Jur\v{c}o, Baas--B\"okstedt--Kro, and others generalizing this result to 
topological 2-groups and even topological 2-categories.  We explain
various viewpoints on topological 2-groups and the \v{C}ech cohomology 
$\check{H}^1(M,\G)$ with coefficients in a topological 2-group $\G$, 
also known as `nonabelian cohomology'.  Then we give an elementary proof 
that under mild conditions on $M$ and $\G$ there is a bijection
$\check{H}^1(M,\G) \cong [M,B|\G|]$ where $B|\G|$ is the classifying 
space of the geometric realization of the nerve of $\G$.  Applying 
this result to the `string 2-group' $\String(G)$ of a simply-connected 
compact simple Lie group $G$, it follows that principal 
$\String(G)$-2-bundles have rational characteristic classes coming from 
elements of $H^*(BG,\mathbb{Q})/\langle c \rangle$, where $c$ is any 
generator of $H^4(BG,\mathbb{Q})$.
  
\end{abstract}

\section{Introduction} 
\label{introduction}

Recent work in higher gauge theory has revealed the importance of
categorifying the theory of bundles and considering `2-bundles', where
the fiber is a topological category instead of a topological space
\cite{BaezSchreiber}.  These structures show up not only in
mathematics, where they form a useful generalization of nonabelian
gerbes \cite{BreenMessing:2001}, but also in physics, where they can
be used to describe parallel transport of strings \cite{SSS,Schreiber}.

The concepts of `\v{C}ech cohomology' and `classifying space' play a
well-known and fundamental role in the theory of bundles.  For 
any topological group $G$, principal $G$-bundles over a space $M$ 
are classified by the first \v{C}ech cohomology of $M$ with 
coefficients in $G$.  Furthermore, under some mild conditions, 
these \v{C}ech cohomology classes are in 1-1 correspondence with 
homotopy classes of maps from $M$ to the classifying space $BG$.  
This lets us define characteristic classes for bundles, coming from
cohomology classes for $BG$.

All these concepts and results can be generalized from bundles to
2-bundles.  Bartels \cite{Bartels:2004} has defined principal
$\G$-2-bundles where $\G$ is a `topological 2-group': roughly
speaking, a categorified version of a topological group.  Furthermore,
his work shows how principal $\G$-2-bundles over $M$ are classified by
$\check{H}^1(M,\G)$, the first \v{C}ech cohomology of $M$ with
coefficients in $\G$.  This form of cohomology, also known as
`nonabelian cohomology', is familiar from work on nonabelian gerbes
\cite{Breen1,Giraud}.

In fact, under mild conditions on $\G$ and $M$, there is a 1-1
correspondence between $\check{H}^1(M,\G)$ and the set of homotopy 
classes of maps from $M$ to a certain space $B|\G|$: the classifying
space of the geometric realization of the nerve of $\G$.   So, 
$B|\G|$ serves as a classifying space for the topological 2-group
$\G$!  This paper seeks to provide an introduction to topological 
2-groups and nonabelian cohomology leading up to a self-contained 
proof of this fact.

In his pioneering work on this subject, Jur\v{c}o \cite{Jurco} 
asserted that a certain space homotopy equivalent to ours is a 
classifying space for the first \v{C}ech cohomology with 
coefficients in $\G$.  However, there are some gaps in his argument for 
this assertion (see Section \ref{remarks} for details).

Later, Baas, B\"okstedt and Kro \cite{BBK} gave the definitive
treatment of classifying spaces for 2-bundles.  For any `good'
topological 2-category $\C$, they construct a classifying space 
$B\C$.  They then show that for any space $M$ with the homotopy
type of a CW complex, concordance classes of `charted $\C$-2-bundles' 
correspond to homotopy classes of maps from $M$ to $B\C$.
In particular, a topological 2-group is just a topological 
2-category with one object and with all morphisms and 2-morphisms 
invertible --- and in this special case, their result \textit{almost} 
reduces to the fact mentioned above.

There are some subtleties, however.  
Most importantly, while their `charted $\C$-2-bundles' 
reduce precisely to our principal $\G$-2-bundles, they classify
these 2-bundles up to concordance, while we classify them up to 
a superficially different equivalence relation.  Two $\G$-2-bundles 
over a space $X$ are `concordant' if they are restrictions of some 
$\G$-2-bundle over $X \times [0,1]$ to the two ends $X \times \{0\}$ 
and $X \times \{1\}$.  This makes it easy to see that homotopic maps 
from $X$ to the classifying space define concordant $\G$-2-bundles.   
We instead consider two $\G$-2-bundles to be equivalent if their 
defining \v{C}ech 1-cocycles are cohomologous.  In this approach, 
some work is required to show that homotopic maps from $X$ to the 
classifying space define equivalent $\G$-2-bundles.  {\em A priori}, 
it is not obvious that two $\G$-2-bundles are equivalent
in this \v{C}ech sense if and only if they are concordant.  However, 
since the classifying space of Baas, B\"okstedt and Kro is homotopy
equivalent to the one we use, it follows from our work that these equivalence 
relations are the same --- at least given $\G$ and $M$ satisfying 
the technical conditions of both their result and ours.

We also discuss an interesting example: the `string 2-group' 
$\String(G)$ of a simply-connected compact simple Lie group $G$ 
\cite{BCSS,Henriques}.  As its name suggests, this 2-group is of 
special interest in physics.  Mathematically, a key fact is that 
$|\String(G)|$ --- the geometric realization of the nerve of 
$\String(G)$ --- is the 3-connected cover of $G$.   Using this, 
one can compute the rational cohomology of $B|\String(G)|$.  
This is nice, because these cohomology classes give `characteristic 
classes' for principal $\G$-2-bundles, and when $M$ is a manifold
one can hope to compute these in terms of a connection and its 
curvature, much as one does for ordinary principal bundles with a 
Lie group as structure group.

Section \ref{summary} is an overview, starting with a review of
the classic results that people are now categorifying.  Section 
\ref{2groups} reviews four viewpoints on topological 2-groups.  
Section \ref{nonabelian} explains nonabelian cohomology with
coefficients in a topological 2-group.  Finally, in Section \ref{proofs} 
we prove the results stated in Section \ref{summary}, and comment
a bit further on the work of Jur\v{c}o and Baas--B\"okstedt--Kro.

\section{Overview}
\label{summary}

Once one knows about `topological 2-groups', it is irresistibly tempting 
to generalize all ones favorite results about topological groups to 
these new entities.  So, let us begin with a quick review of some 
classic results about topological groups and their classifying spaces.

Suppose that $G$ is a topological group.  The \v{C}ech cohomology
$\check{H}^1(M,G)$ of a topological space $M$ with coefficients in $G$
is a set carefully designed to be in 1-1 correspondence with the set
of isomorphism classes of principal $G$-bundles on $M$.  Let us recall
how this works.

First suppose $\cU = \{U_i\}$ is an open cover of $M$ and $P$ is a principal
$G$-bundle over $M$ that is trivial when restricted to each open set
$U_i$.  Then by comparing local trivialisations of $P$ over $U_i$ and
$U_j$ we can define maps $g_{ij}\colon U_i\cap U_j \to G$: the
transition functions of the bundle.  On triple intersections 
$U_i\cap U_j\cap U_k$, these maps satisfy a cocycle condition:
$$ 
g_{ij}(x)g_{jk}(x) = g_{ik}(x) 
$$ 
A collection of maps $g_{ij}\colon U_i\cap U_j \to G$ satisfying 
this condition is called a `\v{C}ech 1-cocycle' subordinate to the
cover $\cU$.  Any such 1-cocycle defines a principal $G$-bundle 
over $M$ that is trivial over each set $U_i$.

Next, suppose we have two principal $G$-bundles over $M$ that are
trivial over each set $U_i$, described by \v{C}ech 1-cocycles $g_{ij}$ and
$g'_{ij}$, respectively.  These bundles are isomorphic if and only if
for some maps $f_i \maps U_i \to G$ we have
$$ 
g_{ij}(x) f_j(x) = f_i(x) g'_{ij}(x) 
$$
on every double intersection $U_i \cap U_j$.  In this case
we say the \v{C}ech 1-cocycles are `cohomologous'.  We
define $\check{H}^1(\cU,G)$ to be the quotient of the set of 
\v{C}ech 1-cocycles subordinate to $\cU$ by this equivalence
relation.  

Recall that a `good' cover of $M$ is an open cover $\cU$ for which all
the non-empty finite intersections of open sets $U_i$ in $\cU$ are
contractible.  We say a space $M$ \textbf{admits good covers} if any
cover of $M$ has a good cover that refines it.  For example, any
(paracompact Hausdorff) smooth manifold admits good covers, as does any
simplicial complex.

If $M$ admits good covers, $\check{H}^1(\cU,G)$ is independent
of the choice of good cover $\cU$.  So, we can denote it simply by 
$\check{H}^1(M,G)$.  Furthermore, this set $\check{H}^1(M,G)$ is in 
1-1 correspondence with the set of isomorphism classes of 
principal $G$-bundles over $M$.  The reason is that we can always
trivialize any principal $G$-bundle over the open sets in a good
cover.  

For more general spaces, we need to define the \v{C}ech cohomology
more carefully.  If $M$ is a paracompact Hausdorff space, we can
define it to be the limit 
\[    \check{H}^1(M,G) = \varinjlim_{\cU}
\check{H}^1(\cU,G)  \]
over all open covers, partially ordered by refinement.  

It is a classic result in topology that $\check{H}^1(M,G)$ can be 
understood using homotopy theory with the help of Milnor's construction
\cite{Dold,Milnor} of the classifying space $BG$: 

\setcounter{theorem}{-1}
\begin{theorem} 
\label{milnor_theorem}
Let $G$ be a topological group.  Then there is a topological space 
$BG$ with the property that for any paracompact Hausdorff space 
$M$, there is a bijection 
$$ 
\check{H}^1(M,G) \cong [M,BG] 
$$ 
\end{theorem}

\noindent
Here $[X,Y]$ denotes the set of homotopy classes of maps from $X$ into
$Y$.  The topological space $BG$ is called the \textbf{classifying space} 
of $G$.  There is a canonical principal $G$-bundle on $BG$, called the
universal $G$-bundle, and the theorem above is usually understood as
the assertion that every principal $G$-bundle $P$ on $M$ is obtained
by pullback from the universal $G$-bundle under a certain map $M\to
BG$ (the classifying map of $P$).

Now let us discuss how to generalize all these results to topological 
2-groups.  First of all, what is a `2-group'?  It is like a group, 
but `categorified'. While a group is a \emph{set}
equipped with \emph{functions} describing multiplication and inverses,
and an identity \emph{element}, a 2-group is a \emph{category}
equipped with \emph{functors} describing multiplication and inverses,
and an identity \emph{object}.  Indeed, 2-groups are also known as
`categorical groups'.

A down-to-earth way to work with 2-groups involves treating them as
`crossed modules'.  A crossed module consists of a pair of groups $H$
and $G$, together with a homomorphism $t\colon H\to G$ and an action
$\alpha$ of $G$ on $H$ satisfying two conditions, equations
\eqref{equivariance} and \eqref{peiffer} below.  Crossed modules were
introduced by J.\ H.\ C.\ Whitehead \cite{Whitehead} without the aid
of category theory.  Mac Lane and Whitehead \cite{MW} later proved
that just as the fundamental group captures all the homotopy-invariant
information about a connected pointed homotopy 1-type, a crossed
module captures all the homotopy-invariant information about a
connected pointed homotopy 2-type.  By the 1960s it was clear to
Verdier and others that crossed modules are essentially the same as
categorical groups.  The first published proof of this may be due to Brown
and Spencer \cite{BS}.

Just as one can define principal $G$-bundles over a space $M$ for any
topological group $G$, one can define `principal $\G$-2-bundles' over
$M$ for any topological 2-group $\G$.  Just as a principal $G$-bundle
has a copy of $G$ as fiber, a principal $\G$-2-bundle has a copy of
$\G$ as fiber.  Readers interested in more details are urged to read
Bartels' thesis, available online \cite{Bartels:2004}.  We shall have
nothing else to say about principal $\G$-2-bundles except that they
are classified by a categorified version of \v{C}ech cohomology,
denoted $\check{H}^1(M,\cG)$.

As before, we can describe this categorified \v{C}ech cohomology as a
set of cocycles modulo an equivalence relation.  Let $\cU$ be a cover
of $M$.  If we think of the $2$-group $\cG$ in terms of its associated
crossed module $(G,H,t,\alpha)$, then a cocycle subordinate to $\cU$
consists (in part) of maps $g_{ij}\colon U_i\cap U_j\to G$ as before.
However, we now `weaken' the cocycle condition and only require that
\begin{equation}
\label{triangle.eq}
t(h_{ijk}) g_{ij}g_{jk} = g_{ik} 
\end{equation}
for some maps $h_{ijk}\colon U_i\cap U_j\cap U_k\to H$.  These maps
are in turn required to satisfy a cocycle condition of their own on
quadruple intersections, namely
\begin{equation}
\label{tetrahedron.eq}
\a(g_{ij})(h_{jkl})h_{ijl} = h_{ijk}h_{ikl}
\end{equation}
where $\alpha$ is the action of $G$ on $H$.  This mildly
intimidating equation will be easier to understand when we draw it as
a commuting tetrahedron --- see equation \eqref{tetrahedron} in 
Section~\ref{nonabelian}.  The pair $(g_{ij},h_{ijk})$ is called a 
$\cG$-valued \textbf{\v{C}ech 1-cocycle} subordinate to $\cU$.  

Similarly, we say two cocycles $(g_{ij},h_{ijk})$
and $(g'_{ij},h'_{ijk})$ are \textbf{cohomologous} if
\begin{equation}
\label{square.eq}
t(k_{ij}) g_{ij} f_j = f_i g'_{ij} 
\end{equation}
for some maps $f_i \maps U_i \to G$ and $k_{ij} \maps U_i \cap U_j \to H$,
which must make a certain prism commute --- see equation
\eqref{prism}.   We define $\check{H}^1(\cU,\G)$ to be the set of 
cohomology classes of $\G$-valued \v{C}ech 1-cocycles.  
To capture the entire cohomology set $\check{H}^1(M,\cG)$, 
we must next take a limit of the sets $\check{H}^1(\cU,\G)$ 
as $\cU$ ranges over all covers of $M$.  For more details we refer 
to Section~\ref{nonabelian}. 

Theorem \ref{milnor_theorem} generalizes nicely from topological 
groups to topological 2-groups.  But, following the usual tradition
in algebraic topology, we shall henceforth work in the category of 
$k$-spaces, i.e., compactly generated weak Hausdorff spaces.  So, 
by `topological space' we shall always mean a $k$-space,
and by `topological group' we shall mean a group object in the
category of $k$-spaces.  
  
\begin{theorem} 
\label{classifying}
Suppose that $\G$ is a well-pointed topological 2-group and
$M$ is a paracompact Hausdorff space admitting good covers.  Then
there is a bijection
$$ \check{H}^1(M,\G) \cong [M,B|\G|]
$$ 
where the topological group $|\G|$ is the geometric realization of
the nerve of $\G$.
\end{theorem} 
\noindent
One term here requires explanation.  A topological group $G$ is 
said to be `well pointed' if $(G,1)$ is an NDR pair, or in other words 
if the inclusion $\{1\}\hookrightarrow G$ is a closed cofibration.  
We say that a topological 2-group $\G$ is {\bf well pointed} if 
the topological groups $G$ and $H$ in its corresponding crossed 
module are well pointed.  
For example, any `Lie 2-group' is well pointed: a topological 2-group is 
called a {\bf Lie 2-group} if $G$ and $H$ are Lie groups and the maps 
$t,\alpha$ are smooth.
More generally, any `Fr\'echet Lie 2-group' \cite{BCSS} is well pointed.  
We explain the importance of this notion in Section
\ref{proof.thm.1}.

Bartels \cite{Bartels:2004} has already considered two examples of principal
$\G$-2-bundles, corresponding to abelian gerbes and nonabelian gerbes.
Let us discuss the classification of these before turning to a third,
more novel example.

For an abelian gerbe \cite{Brylinski}, we first choose an abelian
topological group $H$ --- in practice, usually just $\U(1)$.  Then, we
form the crossed module with $G = 1$ and this choice of $H$, with $t$
and $\alpha$ trivial.  The corresponding topological 2-group deserves
to be called $H[1]$, since it is a `shifted version' of $H$.  Bartels
shows that the classification of abelian $H$-gerbes matches the
classification of $H[1]$-2-bundles.  It is well-known that
$$  |H[1]| \cong BH   $$
so the classifying space for abelian $H$-gerbes is 
$$ B|H[1]| \cong B(BH)  $$
In the case $H = \U(1)$, this classifying space is just $K(\ZZ,3)$.  
So, in this case, we recover the well-known fact that abelian 
$\U(1)$-gerbes over $M$ are classified by
\[          [M, K(\ZZ,3)] \cong H^3(M,\ZZ) \]
just as principal $\U(1)$ bundles are classified by $H^2(M,\ZZ)$.

For a nonabelian gerbe \cite{Breen1,GS,Giraud}, we fix any topological
group $H$.  Then we form the crossed module with $G = \Aut(H)$ and
this choice of $H$, where $t \maps H \to G$ sends each element of $H$
to the corresponding inner automorphism, and the action of $G$ on $H$
is the tautologous one.  This gives a topological 2-group called
$\AUT(H)$.  Bartels shows that the classification of nonabelian
$H$-gerbes matches the classification of $\AUT(H)$-2-bundles.  It
follows that, under suitable conditions on $H$, nonabelian $H$-gerbes are classified by homotopy classes
of maps into $B|\AUT(H)|$.

A third application of Theorem \ref{classifying} arises when
$G$ is a simply-connected compact simple Lie group.  For any such
group there is an isomorphism $H^3(G,\ZZ) \cong \ZZ$ and the generator
$\nu\in H^3(G,\ZZ)$ transgresses to a characteristic class $c\in
H^4(BG,\ZZ) \cong \ZZ$.  Associated to $\nu$ is a map $G\to
K(\ZZ,3)$ and it can be shown that the homotopy fiber of this
can be given the structure of a topological group $\hat{G}$.  
This group $\hat{G}$ is the 3-connected cover of $G$.  
When $G = \Spin(n)$, this group $\hat{G}$ is known as
$\String(n)$.  In general, we might call $\hat{G}$ the {\bf string
group} of $G$.  Note that until one picks a specific construction for
the homotopy fiber, $\hat{G}$ is only defined up to homotopy --- or
more precisely, up to equivalence of $A_\infty$-spaces.

In \cite{BCSS}, under the above hypotheses on $G$, a topological
2-group subsequently dubbed the \textbf{string 2-group} of $G$ was
introduced.  Let us denote this by $\String(G)$.  A key
result about $\String(G)$ is that the topological group $|\String(G)|$
is equivalent to $\hat{G}$.  By construction $\String(G)$ is a
Fr\'echet Lie 2-group, hence well pointed.  So, from Theorem
\ref{classifying} we immediately conclude:

\begin{corollary}
\label{equivalence}
Suppose that $G$ is a simply-connected compact simple Lie group.
Suppose $M$ is a paracompact Hausdorff space admitting good
covers.  Then there are bijections between the following sets:
\begin{itemize}
\item
the set of equivalence classes of principal $\String(G)$-2-bundles over $M$,
\item 
the set of isomorphism classes of principal $\hat{G}$-bundles over $M$,
\item
$\check{H}^1(M,\String(G))$,
\item
$\check{H}^1(M,\hat{G})$,
\item
$[M, B\hat{G}]$.
\end{itemize}
\end{corollary}

One can describe the rational cohomology of $B\hat{G}$ in terms of the
rational cohomology of $BG$, which is well-understood.  The following
result was pointed out to us by Matt Ando \cite{Ando}, and later 
discussed by Greg Ginot \cite{Ginot}:

\begin{theorem} 
\label{string}
Suppose that $G$ is a simply-connected compact simple Lie group, and
let $\hat{G}$ be the string group of $G$.  Let $c\in H^4(BG,\QQ)
= \QQ$ denote the transgression of the generator $\nu\in H^3(G,\QQ) =
\QQ$.  Then there is a ring isomorphism
$$ 
H^*(B\hat{G},\QQ) \cong H^*(BG,\QQ)/\langle c\rangle  
$$ 
where $\langle c \rangle$ is the ideal generated by $c$.
\end{theorem} 

As a result, we obtain characteristic classes for $\String(G)$-2-bundles:

\begin{corollary} 
\label{characteristic}
Suppose that $G$ is a simply-connected compact simple Lie group and
$M$ is a paracompact Hausdorff space admitting good covers.  Then an
equivalence class of principal $\String(G)$-2-bundles over $M$ 
determines a ring homomorphism
$$ H^*(BG,\QQ)/\langle c \rangle \to H^*(M,\QQ) $$
\end{corollary} 

To see this, we use Corollary \ref{equivalence} to 
reinterpret an equivalence class of principal $\G$-2-bundles over $M$ 
as a homotopy class of maps $f \maps M \to B|\G|$.  Picking any 
representative $f$, we obtain a ring homomorphism
$$ f^\ast \maps H^*(B|\G|,\QQ) \to H^*(M,\QQ). $$
This is independent of the choice of representative.  Then,
we use Theorem \ref{string}.

It is a nice problem to compute the rational characteristic classes 
of a principal $\String(G)$-2-bundle over a manifold using de Rham
cohomology.  It should be possible to do this using the curvature of
an arbitrary connection on the 2-bundle, just as for ordinary principal
bundles with a Lie group as structure group.  Sati, Schreiber and
Stasheff \cite{SSS} have recently made excellent progress on solving
this problem and its generalizations to $n$-bundles for higher $n$.  

\section{Topological $2$-Groups} 
\label{2groups}

In this section we recall four useful perspectives on topological
$2$-groups.  For a more detailed account, we refer the reader to
\cite{HDA5}.  

Recall that for us, a `topological space'
really means a $k$-space, and a `topological group' really means a
group object in the category of $k$-spaces.
A \textbf{topological $2$-group} is a groupoid in the category of
topological groups.  In other words, it is a groupoid $\G$ where
the set $\Ob(\G)$ of objects and the set $\Mor(\G)$ of morphisms 
are each equipped with the structure of a topological group such that 
the source and target maps $s,t\colon \Mor(\G)\to \Ob(\G)$, the 
map $i\colon \Ob(\G)\to \Mor(\G)$ assigning each object its identity
morphism, the composition map 
$\circ\colon \Mor(\G)\times_{\Ob(\G)}\Mor(\G)\to \Mor(\G)$,
and the map sending each morphism to its inverse are all continuous 
group homomorphisms.   

Equivalently, we can think of a topological 2-group as a group in
the category of topological groupoids.  A \textbf{topological groupoid}
is a groupoid $\G$ where $\Ob(\G)$ and $\Mor(\G)$ are topological
spaces (or more precisely, $k$-spaces) and all the groupoid
operations just listed are continuous maps.  
We say that a functor $f\colon \G\to
\G'$ between topological groupoids is \textbf{continuous} if
the maps $f\colon \Ob(\G)\to \Ob(\G')$ and $f\colon \Mor(\G)\to
\Mor(\G')$ are continous.
A group in
the category of topological groupoids is such a thing equipped with
continuous functors $m \maps \G \times \G \to \G$, $\mathrm{inv} \maps \G 
\to \G$ and a unit object $1 \in \G$ satisfying the usual group axioms,
written out as commutative diagrams.  

This second viewpoint is useful because
any topological groupoid $\G$ has a `nerve' $N\G$, a simplicial 
space where the space of $n$-simplices consists of composable strings of
morphisms
$$  
x_0 \stackrel{f_1}{\longrightarrow} x_1 \stackrel{f_2}{\longrightarrow}
\cdots \stackrel{f_{n-1}}{\longrightarrow} 
x_{n-1} \stackrel{f_n}{\longrightarrow} x_n 
$$
Taking the geometric realization of this nerve, we obtain a 
topological space which we denote as $|\G|$ for short.
If $\G$ is a topological 2-group, its nerve inherits a group
structure, so that $N\G$ is a topological simplicial group.  This in
turn makes $|\G|$ into a topological group.  

A third way to understand topological 2-groups is to view them as
topological crossed modules.  Recall that a {\bf topological crossed module} 
$(G,H,t,\a)$ consists of topological groups $G$ and $H$ together 
with a continuous homomorphism 
$$ 
t\colon H\to G\qquad 
$$ 
and a continuous action
$$
\begin{array}{rcl}
\alpha \maps G \times H & \to & H \\
   (g,h) & \mapsto & \alpha(g) h
\end{array}
$$
of $G$ as automorphisms of $H$,
satisfying the following two identities: 
\begin{align} 
\label{equivariance}
& t(\a(g)(h)) = gt(h)g^{-1} \\ 
\label{peiffer}
& \a(t(h))(h') = hh'h^{-1}. 
\end{align} 
The first equation above implies that the map $t\colon H\to G$ is equivariant 
for the action of $G$ on $H$ defined by $\a$ 
and the action of $G$ on itself by conjugation.  The second equation
is called the \textbf{Peiffer identity}.  When no confusion is likely to 
result, we will sometimes denote the 2-group corresponding
to a crossed module $(G,H,t,\a)$ simply by $H \to G$.  

Every topological crossed module determines a topological 
$2$-group and vice versa.  Since there are some choices of
convention involved in this construction, we briefly review it 
to fix our conventions.   Given a topological crossed module
$(G,H,t,\a)$, we define a topological 2-group $\G$ as follows.
First, define the group $\Ob(\G)$ of objects of $\G$ and the group 
$\Mor(\G)$ of morphisms of $\G$ by
$$ 
\Ob(\G) = G,\qquad \Mor(\G) = H \rtimes G
$$ 
where the semidirect product $H\rtimes G$ is formed using the left 
action of $G$ on $H$ via $\a$: 
$$ 
(h,g)\cdot (h',g') = (h\a(g)(h'), gg') 
$$ 
for $g,g'\in G$ and $h,h'\in H$.  The source and target of 
a morphism $(h,g) \in \Mor(\G)$ are defined by
$$ 
s(h,g) = g\qquad \text{and}\qquad t(h,g) = t(h)g 
$$ 
(Denoting both the target map $t \maps
\Mor(\G) \to \Ob(\G)$ and the homomorphism $t\maps H\to G$ by
the same letter should not cause any problems, since the first
is the restriction of the second to $H \subseteq \Mor(\G)$.)
The identity morphism of an object $g \in \Ob(G)$ is defined by
$$i(g) = (1,g). $$
Finally, the composite of the morphisms 
$$\alpha = (h,g) \maps g \to t(h)g \qquad \text{and}\qquad
\beta = (h',t(h)g) \maps t(h)g \to t(h'h)g'$$ 
is defined to be
$$ \beta \circ \alpha =  (h'h, g) \maps g \to t(h'h) g $$
It is easy to check that with these definitions, $\G$ is a $2$-group.  
Conversely, given a topological $2$-group $\G$, we define a crossed module 
$(G,H,t,\a)$ by setting $G$ to be $\Ob(\G)$, $H$ to be 
$\ker(s)\subset \Mor(\G)$, 
$t$ to be the restriction of the target homomorphism 
$t\colon \Mor(\G)\to \Ob(\G)$ to the subgroup $H\subset \Mor(\G)$, and
setting
$$\a(g)(h) = i(g) h i(g)^{-1}  $$

If $G$ is any topological group then there is a topological
crossed module $1\to G$ where $t$ and $\alpha$ are trivial.
The underlying groupoid of the corresponding topological 2-group 
has $G$ as its space of objects, and only identity morphisms.  We
sometimes call this 2-group the \textbf{discrete} topological 2-group
associated to $G$ --- where `discrete' is used in the sense of
category theory, not topology!  

At the other extreme, if $H$ is a topological group then it follows
from the Peiffer identity that $H\to 1$ can be made into topological
crossed module if and only if $H$ is abelian, and then in a unique
way.  This is because a groupoid with one object and $H$ as
morphisms can be made into a 2-group precisely when $H$ is abelian.
We already mentioned this 2-group in the previous section, where we
called it $H[1]$.

We will also need to talk about homomorphisms of $2$-groups.  
We shall understand these in the strictest possible sense.  So, 
we say a \textbf{homomorphism} of topological $2$-groups is a functor 
such that $f \maps \Ob(\G)\to \Ob(\G')$ and $f\maps \Mor(\G)\to 
\Mor(\G')$ are both continuous homomorphisms of topological groups.  We can also describe 
$f$ in terms of the crossed modules $(G,H,t,\a)$ and $(G',H',t',\a')$ 
associated to $\G$ and $\G'$ respectively.  In these terms the data 
of the functor $f$ is described by the commutative diagram 
$$ 
\xymatrix{ 
H \ar[r]^-{f} \ar[d]_-t & H' \ar[d]^-{t'} \\ 
G \ar[r]_-{f} & G' } 
$$ 
where the upper $f$ denotes the restriction of $f \maps \Mor(\G) 
\to \Mor(\G')$ to a map from $H$ to $H'$.  (We are using 
$f$ to mean several different things, but this makes the notation
less cluttered, and should not cause any confusion.)
The maps $f \maps G \to G'$ and $f \maps H \to H'$
must both be continuous homomorphisms, and moreover must
satisfy an equivariance property with respect to the actions of $G$ on 
$H$ and $G'$ on $H'$: we have 
$$ 
f(\a(g)(h)) = \a(f(g))(f(h)) 
$$ 
for all $g\in G$ and $h\in H$.  

Finally, we will need to talk about short exact sequences of topological
groups and $2$-groups.   Here the topology is important.  
If $G$ is a topological group and $H$ is a normal
topological subgroup of $G$, then we can define an action of $H$ on $G
$ by right translation.  In some circumstances, the projection 
$G\to G/H$ is a Hurewicz fibration.
For instance, this is the case if $G$ is a Lie group and $H$ is a closed 
normal subgroup of $G$.  We define a 
{\bf short exact sequence} of topological groups to be a sequence 
$$1 \to H \to G \to K \to 1$$ 
of topological groups and continuous homomorphisms 
such that the underlying sequence of groups is exact and 
the map underlying the homomorphism $G\to K$ is a Hurewicz fibration.

Similarly, we define a {\bf short exact sequence} of topological
2-groups to be a sequence
$$1\to \G'\to \G\to \G''\to 1$$
of topological 2-groups and continuous homomorphisms between them
such that both the resulting sequences
\[
1\to \Ob(\G')\to \Ob(\G)\to \Ob(\G'')\to 1 
\]
\[
1\to \Mor(\G')\to \Mor(\G)\to \Mor(\G'')\to 1
\]
are short exact sequences of topological groups.
Again, we can interpret this in terms of the associated crossed modules: 
if $(G,H,t,\a)$, $(G',H',t',\a')$ and $(G'',H'',t'',\a'')$ denote the 
associated crossed modules, then it can be shown that the sequence of topological 2-groups
$1\to \G'\to \G\to \G''\to 1$
is exact if and only if both rows in the commutative diagram 
$$ 
\xymatrix{ 
1\ar[r] & H' \ar[r] \ar[d] & H\ar[r] \ar[d] & H'' \ar[d] \ar[r] & 1 \\ 
1\ar[r] & G' \ar[r] & G\ar[r] & G'' \ar[r] &  1} 
$$ 
are short exact sequences of topological groups.  In this situation we
also say we have a short exact sequence of topological crossed
modules.  

At times we shall also need a fourth viewpoint on topological
2-groups: they are strict topological 2-groupoids with a single
object, say $\bullet$.  In this approach, what we had been calling `objects'
are renamed `morphisms', and what we had been calling `morphisms'
are renamed `2-morphisms'.  This verbal shift can be confusing, so we will 
not engage in it!  However, the 2-groupoid viewpoint is very handy
for diagrammatic reasoning in nonabelian cohomology.  
We draw $g \in \Ob(\G)$ as an arrow: 
\[
\xymatrix{
   \bullet \ar[rr]^{g}
&& \bullet
}
\]
and draw $(h,g) \in \Mor(\G)$ as a bigon:
\[
\xymatrix{
   \bullet \ar@/^1pc/[rr]^{g}_{}="0"
           \ar@/_1pc/[rr]_{g'}_{}="1"
           \ar@{=>}"0";"1"^{h}
&& \bullet
}
\]
where $g'$ is the target of $(h,g)$, namely $t(h) g$.
With our conventions, horizontal composition of
2-morphisms is then given by:
\[ 
\xymatrix{ 
\bullet \ar@/^1pc/[rr]^-{g_1}_{}="0" 
            \ar@/_1pc/[rr]_{g'_1}_{}="1" 
            \ar@{=>}"0";"1"^{h_1} 
& & \bullet 
           \ar@/^1pc/[rr]^-{g_2}_{}="2" 
            \ar@/_1pc/[rr]_{g'_2}_{}="3" 
            \ar@{=>}"2";"3"^{h_2} 
& & \bullet 
}
\quad = \quad 
\xymatrix{
\bullet \ar@/^2pc/[rrrr]^-{g_1g_2}_{}="0" 
            \ar@/_2pc/[rrrr]_{g'_1g'_2}_{}="1" 
            \ar@{=>}"0";"1"^{h_1\a(g_1)(h_2)}
& & & & \bullet            
}
\] 
while vertical composition is given by:
\[ 
\xymatrix{ 
\bullet \ar@/^2pc/[rr]^-{g}_{}="0" 
            \ar[rr]_{}="1" 
            \ar@/_2pc/[rr]_{g'}_{}="2" 
            \ar@{=>} "0";"1"^{h} 
            \ar@{=>} "1";"2"^{h'} 
& & \bullet 
} 
\quad = \quad 
\xymatrix{
\bullet \ar@/^1pc/[rr]^-{g}_{}="0" 
            \ar@/_1pc/[rr]_{g'}_{}="1" 
            \ar@{=>}"0";"1"^{h'h}
& & \bullet            
}
\]

\section{Nonabelian Cohomology} 
\label{nonabelian}

In Section \ref{summary} we gave a quick sketch of nonabelian
cohomology.  The subject deserves a more thorough and more conceptual
explanation.  

As a warmup, consider the \v{C}ech cohomology of a space $M$ with
coefficients in a topological group $G$.  In this case, Segal
\cite{Segal1} realized that we can reinterpret a \v{C}ech 1-cocycle as
a \textit{functor}.  Suppose $\cU$ is an open cover of $M$.  Then there
is a topological groupoid $\hat{\cU}$ whose objects are pairs $(x,i)$ with
$x\in U_i$, and with a single morphism from $(x,i)$ to
$(x,j)$ when $x\in U_i\cap U_j$, and none otherwise.  
We can also think of $G$ as a topological groupoid with a single
object $\bullet$.  Segal's key observation was that a continuous functor
$$        g \maps \hat{\cU} \to G    $$
is the same as a normalized \v{C}ech 1-cocycle subordinate to $\cU$.

To see this, note that a functor $g \maps \hat{\cU} \to G$ maps each
object of $\hat{\cU}$ to $\bullet$, and each morphism $(x,i)\to (x,j)$
to some $g_{ij}(x) \in G$.  For the functor to preserve composition,
it is necessary and sufficient to have the cocycle equation
$$ 
g_{ij}(x) g_{jk}(x) = g_{ik}(x) 
$$ 
We can draw this suggestively as a commuting triangle in the groupoid
$G$:
\[
\xy 
(0,0)*+{\bullet}="1"; 
(-12,-20.78)*+{\bullet}="2"; 
(12,-20.78)*+{\bullet}="3"; 
{\ar^-{g_{ij}(x)} "2";"1"}; 
{\ar^-{g_{jk}(x)} "1";"3"}; 
{\ar_-{g_{ik}(x)} "2";"3"}; 
\endxy
\]
For the functor to preserve identities, it is necessary and sufficient
to have the normalization condition $g_{ii}(x) = 1$.

In fact, even more is true: two cocycles $g_{ij}$ and $g_{ij}'$
subordinate to $\cU$ are cohomologous if and only if the corresponding
functors $g$ and $g'$ from $\hat{\cU}$ to $G$ have a continuous
natural isomorphism between them.  To see this, note that $g_{ij}$ and
$g'_{ij}$ are cohomologous precisely when there are maps $f_i \maps
U_i \to G$ satisfying
$$g_{ij}(x) f_j(x) = f_i(x) g'_{ij}(x)$$
We can draw this equation as a commuting square in the groupoid $G$:
\[ 
\xymatrix{ 
\bullet \ar[r]^-{g_{ij}(x)} \ar[d]_-{f_i(x)} & \bullet \ar[d]^-{f_j(x)} \\ 
\bullet \ar[r]_-{g_{ij}'(x)} & \bullet }
\] 
This is precisely the naturality square for a natural isomorphism
between the functors $g$ and $g'$.

One can obtain \v{C}ech cohomology with coefficients in a 2-group by
categorifying Segal's ideas.  Suppose $\G$ is a topological 2-group
and let $(G,H,t,\alpha)$ be the corresponding topological crossed
module.  Now $\G$ is the same as a topological 2-groupoid with one
object $\bullet$.  So, it is no longer appropriate to consider mere
\emph{functors} from $\hat{\cU}$ into $\G$.  Instead, we should
consider \emph{weak $2$-functors}, also known as `pseudofunctors'
\cite{KS}.  For this, we should think of $\hat{\cU}$ as a topological
2-groupoid with only identity 2-morphisms.

Let us sketch how this works.  A weak 2-functor $g \maps \hat{\cU} \to
\G$ sends each object of $\hat{\cU}$ to $\bullet$, and each 1-morphism
$(x,i)\to (x,j)$ to some $g_{ij}(x) \in G$.  However, composition of
1-morphisms is only weakly preserved.  This means the above triangle
will now commute only up to isomorphism:
\[
\xy 
(0,0)*+{\bullet}="1"; 
(-12,-20.78)*+{\bullet}="2"; 
(12,-20.78)*+{\bullet}="3"; 
{\ar^-{g_{ij}} "2";"1"}; 
{\ar^-{g_{jk}} "1";"3"}; 
{\ar_-{g_{ik}} "2";"3"}; 
{\ar@2{->}^-{h_{ijk}} (0,-7)*{};(0,-18)*{}};
\endxy
\]
where for readability we have omitted the dependence on $x\in U_i\cap
U_j\cap U_k$.  Translated into equations, this triangle says that we have
continuous maps $h_{ijk} \maps U_i \cap U_j \cap U_k \to H$ satisfying
\[  
g_{ik}(x) = t(h_{ijk}(x))g_{ij}(x)g_{jk}(x) 
\]
This is precisely equation \eqref{triangle.eq} from Section \ref{summary}.   

For a weak 2-functor, it is not merely true that composition is preserved
up to isomorphism: this isomorphism is also subject to a coherence law.
Namely, the following tetrahedron must commute:  
\begin{equation}
\label{tetrahedron}
\begin{picture}(110,130)
  \includegraphics[scale = 0.85]{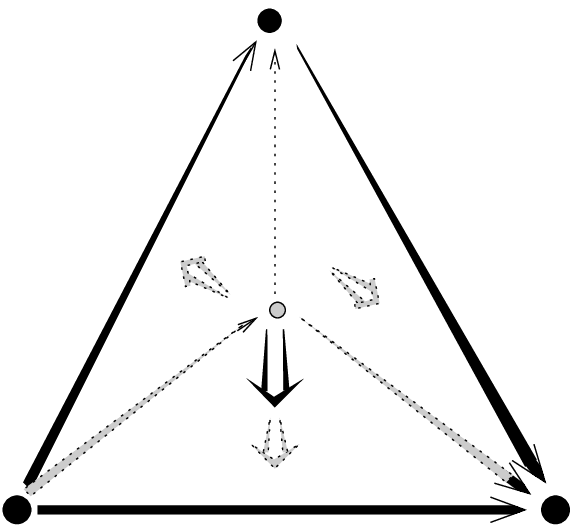}
  \put(-72,90){${}_{g_{jk}}$}
  \put(-40,35){${}_{g_{jl}}$}
  \put(-750,-5){$g_{il}$}
  \put(-125,60){$g_{ik}$}
  \put(-32.5,60){$g_{kl}$}
  \put(-100,35){${}_{g_{ij}}$}
  \put(-65,16){$h_{ijl}$}
  \put(-66.5,40){$h_{ikl}$}
  \put(-60,70){${}_{h_{jkl}}$}
  \put(-100,72.5){${}_{h_{ijk}}$}
  \put(-75,-5){$g_{il}$}
\end{picture}
\end{equation}
where again we have omitted the dependence on $x$.  The commutativity
of this tetrahedron is equivalent to the following equation:
$$
\a(g_{ij})(h_{jkl})h_{ijl} = h_{ijk}h_{ikl}
$$
holding for all $x\in U_i\cap U_j\cap U_k\cap U_l$.  This is
equation \eqref{tetrahedron.eq}.

A weak 2-functor may also preserve identity 1-morphisms only up to isomorphism.
However, it turns out \cite{Bartels:2004} that without loss of generality
we may assume that $g$ preserves identity 1-morphisms \emph{strictly}.  
Thus we have $g_{ii}(x) = 1$ for all $x\in U_i$.  We may also assume
$h_{ijk}(x) = 1$ whenever two or more of the indices $i$, $j$ and $k$
are equal.  Finally, just as for the case of an ordinary topological
group, we require that $g$ is a \textit{continuous} weak 2-functor
We shall not spell this out in detail; suffice it to say that the maps
$g_{ij}\colon U_i\cap U_j\to G$ and $h_{ijk}\colon U_i\cap U_j\cap
U_k\to H$ should be continuous.   We say such continuous
weak 2-functors $g\colon \hat{\cU}\to \G$ are \textbf{\v{C}ech 
1-cocycles} valued in $\G$, subordinate to the cover $\cU$.

We now need to understand when two such cocycles should be considered 
equivalent.  In the case of cohomology with coefficients in an ordinary 
topological group, we saw that two cocycles were cohomologous
precisely when there was a continuous natural isomorphism between 
the corresponding functors.  In our categorified setting we should 
instead use a `weak natural isomorphism', also called a
pseudonatural isomorphism \cite{KS}.  So, we declare two cocycles to be 
{\bf cohomologous} if there is a continuous weak natural isomorphism 
$f\maps g \Rightarrow g'$ between the corresponding weak 2-functors 
$g$ and $g'$.  

In a weak natural isomorphism, the usual naturality square commutes
only up to isomorphism.  So, $f \maps g \Rightarrow g'$ not only sends 
every object $(x,i)$ of $\hat{\cU}$ to some $f_i(x) \in G$, but also
sends every morphism $(x,i)\to(x,j)$ to some $k_{ij}(x) \in H$ filling
in this square:
\[ 
\xy
(-7.5,7.5)*+{\bullet}="1"; 
(-7.5,-7.5)*+{\bullet}="2"; 
(7.5,7.5)*+{\bullet}="3"; 
(7.5,-7.5)*+{\bullet}="4";
{\ar_-{f_{i}} "1";"2"}; 
{\ar^-{g_{ij}} "1";"3"}; 
{\ar^-{f_{j}} "3";"4"};
{\ar_-{g'_{ij}} "2";"4"}; 
{\ar@2{->}^-{k_{ij}} (2,2)*{};(-2,-2)*{}};
\endxy 
\]
Translated into equations, this square says that
\[ 
t(k_{ij}) g_{ij} f_j = f_i g'_{ij}
\]   
This is equation \eqref{square.eq}.

There is also a coherence law that the $k_{ij}$ must satisfy: they
must make the following prism commute:
\begin{eqnarray}
\label{prism}
\begin{picture}(200,270)
\includegraphics{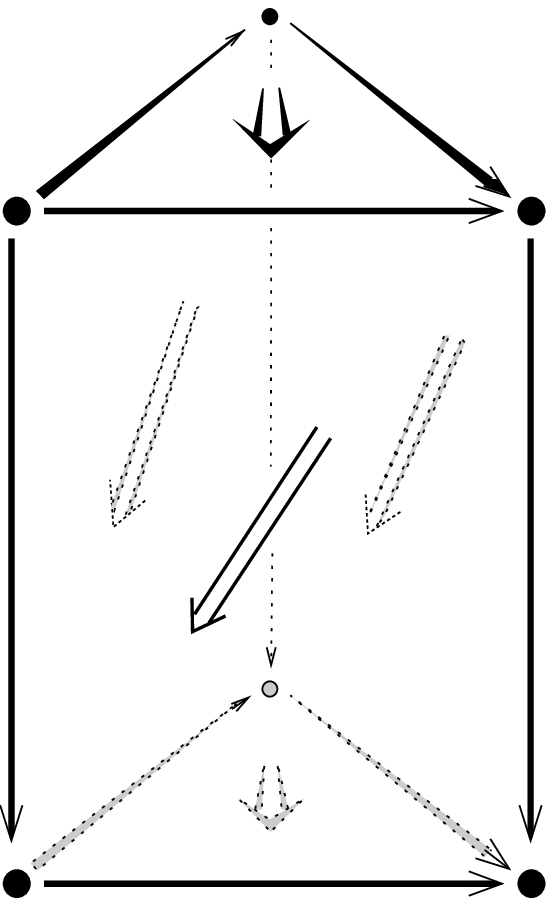}
  \put(-1,100){$f_k$}
  \put(-170,100){$f_i$}
  \put(-80,160){${}_{f_j}$}
  \put(-90,-7){$g'_{ik}$}
  \put(-100,187){$g_{ik}$}
  \put(-128,43){${}_{g'_{ij}}$}
  \put(-128,237){${}_{g_{ij}}$}
  \put(-47,43){${}_{g'_{jk}}$}
  \put(-47,237){${}_{g_{jk}}$}
  \put(-95,16){${}_{h'_{ijk}}$}
  \put(-95,212){${}_{h_{ijk}}$}
  \put(-95,115){$k_{ik}$}
  \put(-132,142){${}_{k_{ij}}$}
  \put(-35,121){${}_{k_{jk}}$}
\end{picture}
\end{eqnarray}
At this point, translating the diagrams
into equations becomes tiresome and unenlightening.
 
It can be shown that this notion of `cohomologousness' of \v{C}ech
1-cocycles $g\maps \hat{\cU}\to \G$ is an equivalence relation.  We
denote by $\check{H}^1(\cU,\G)$ the set of equivalence classes of
cocycles obtained in this way.  In other words, we let
$\check{H}^1(\cU,\G)$ be the set of continuous weak natural
isomorphism classes of continuous weak 2-functors $g \maps \hat{\cU}\to \G$.

Finally, to define $\hat{H}^1(M,\G)$, we need to take all covers into
account as follows.  The set of all open covers of $M$ is a directed
set, partially ordered by refinement.  By restricting cocycles defined
relative to $\cU$ to any finer cover $\cV$, we obtain a map
$\check{H}^1(\cU,\G) \to \check{H}^1(\cV,\G)$.  This allows us 
to define the \v{C}ech cohomology $\check{H}^1(M,\G)$ as a limit:

\begin{definition} Given a topological space $M$ and a topological
2-group $\G$, we define the \textbf{first \v{C}ech cohomology
of} $M$ \textbf{with coefficients in} $\G$ to be
$$ 
\check{H}^1(M,\G) = \varinjlim_{\cU}\check{H}^1(\cU,\G)
$$
\end{definition}
\noindent 
When we want to emphasize the crossed module, we will sometimes use
the notation $\check{H}^1(M,H\to G)$ instead of $\check{H}^1(M,\G)$.
Note that $\check{H}^1(M,\G)$ is a pointed set, pointed by the trivial
cocycle defined relative to any open cover $\{U_i\}$ by $g_{ij} = 1$,
$h_{ijk} = 1$ for all indices $i$, $j$ and $k$.

In Theorem \ref{classifying} we assume $M$ admits good covers, so that
every cover $\cU$ of $M$ has a refinement by a good cover $\cV$.  In
other words, the directed set of good covers of $M$ is cofinal in the
set of all covers of $M$.  As a result, in computing the limit
above, it is sufficient to only consider \emph{good} covers $\cU$.

Finally, we remark that there is a more refined version of the set
$\check{H}^1(M,\G)$ defined using the notion of `hypercover'
\cite{Breen1, Brown, Jardine}.  For a paracompact space $M$ this
refined cohomology set $H^1(M,G)$ is isomorphic to the set
$\check{H}^1(M,G)$ defined in terms of \v{C}ech covers.  While the
technology of hypercovers is certainly useful, and can simplify some
proofs, our approach is sufficient for the applications we have in
mind (see also the remark following the proof of Lemma~\ref{lem2} in
sub-section~\ref{proof.lem.2}).

\section{Proofs}
\label{proofs}

\subsection{Proof of Theorem \ref{classifying}} 
\label{proof.thm.1}

First, we need to distinguish between Milnor's \cite{Milnor}
original construction of a classifying space for a topological group and 
a later construction introduced by Milgram, Segal and Steenrod
\cite{Milgram,Segal1,Steenrod} and further studied by May 
\cite{May1}.  Milnor's construction is very powerful, 
as witnessed by the generality of Theorem \ref{milnor_theorem}.  The
later construction is conceptually more beautiful: for any topological 
group $G$, it constructs $BG$ as the geometric realization of the nerve 
of the topological groupoid with one object associated to $G$.  But, 
here we are performing this construction in the category of $k$-spaces,
rather than the traditional category of topological spaces. 
It also seems to give a slightly weaker result: to obtain a 
bijection 
$$\check{H}^1(M,G) \cong [M,BG]$$ 
all of the above cited works require some extra hypotheses on $G$: 
Segal \cite{Segal2} requires that $G$ be locally contractible; May, 
Milgram and Steenrod require that $G$ be well pointed.  
This extra hypothesis on $G$ is required in the construction of 
the universal principal $G$-bundle $EG$ over $BG$; to ensure 
that the bundle is locally trivial we must make one of the above 
assumptions on $G$.  May's work goes further in this regard: he 
proves that if $G$ is well pointed then $EG$ is a \emph{numerable} principal 
$G$-bundle over $BG$, and hence $EG\to BG$ is a Hurewicz fibration.  

Another feature of this later construction is that $EG$ comes equipped
with the structure of a topological group.  In the work of May and Segal, 
this arises from the fact that $EG$ is the geometric realization of the 
nerve of a topological 2-group.  We need the group structure on $EG$, 
so we will use this later construction rather than Milnor's.  For further 
comparison of the constructions see tom Dieck \cite{tomDieck}.

We prove Theorem \ref{classifying} using three lemmas that are of 
some interest in their own right.  The second, as far as we know, 
is due to Larry Breen:

\begin{lemma} \label{lem1} 
Let $\G$ be any well-pointed topological 2-group, and let
$(G,H,t,\alpha)$ be the corresponding topological crossed module.  
Then:
\begin{enumerate} 
\item $|\G|$ is a well-pointed topological group.
\item There is a topological 
2-group $\hat{\G}$ such that $|\hat{\G}|$ fits into a short
exact sequence of topological groups
$$ 
1\to H \to |\hat{\G}| \stackrel{p}{\to} |\G| \to 1
$$ 
\item $G$ acts continuously via automorphisms on the topological group $EH$,
and there is an isomorphism $|\hat{\G}| \cong G\ltimes EH$.  This exhibits   
$|\G|$ as $G \ltimes_H EH$, the quotient of 
$G\ltimes EH$ by the normal subgroup $H$. 
\end{enumerate} 
\end{lemma}

\begin{lemma} \label{lem2}
If 
$$1\to H \stackrel{t}{\to} G\stackrel{p}{\to} K\to 1$$
is a short exact sequence of topological groups, there is a bijection 
$$ 
\check{H}^1(M,H\to G)\cong \check{H}^1(M,K) 
$$ 
Here $H \to G$ is our shorthand for the 2-group corresponding to the
crossed module $(G,H,t,\alpha)$ where $t$ is the inclusion of the
normal subgroup $H$ in $G$ and $\alpha$ is the action of $G$ by 
conjugation on $H$.
\end{lemma}

\begin{lemma} \label{lem3}
If 
$$ 1 \to \G_0 \stackrel{f}{\to} \G_1 \stackrel{p}{\to} \G_2 \to 1$$ 
is a short exact sequence of topological 2-groups, then 
$$ 
\check{H}^1(M,\G_0) \stackrel{f_\ast}{\to} \check{H}^1(M,\G_1) 
\stackrel{p_\ast}{\to} \check{H}^1(M,\G_2) 
$$ 
is an exact sequence of pointed sets.
\end{lemma}

Given these lemmas the proof of Theorem \ref{classifying} goes as follows.  
Assume that $\G$ is a well-pointed topological 2-group.  
From Lemma~\ref{lem1} we see that $|\G|$ is a well-pointed topological group.  
It follows that we have a bijection
$$\check{H}^1(M,|\G|) \cong [M,B|\G|]$$ 
So, to prove the theorem, it suffices to construct a bijection
$$ 
\check{H}^1(M,\G) \cong \check{H}^1(M,|\G|) 
$$ 
By Lemma~\ref{lem1}, $|\G|$ fits into a short exact sequence 
of topological groups:
$$1\to H\to G\ltimes EH\to |\G|\to 1$$
We can use Lemma~\ref{lem2} to conclude that there is a bijection
$$ 
\check{H}^1(M,H\to G\ltimes EH)
\cong \check{H}^1(M,|\G|)
$$ 
To complete the proof it thus suffices to construct a bijection
$$
\check{H}^1(M,H\to G\ltimes EH) \cong
\check{H}^1(M,\G) 
$$
For this, observe that we have a short exact sequence of topological 
crossed modules:
\[
\xymatrix{ 
1 \ar[r] & 1 \ar[r] \ar[d] & H \ar[d] \ar[r]^-1 & H \ar[d]^-{t} 
\ar[r] & 1 \\ 
1 \ar[r] & EH \ar[r] & G \ltimes EH \ar[r] & G \ar[r] & 1 } 
\]
So, by Lemma~\ref{lem3}, we have an exact sequence of sets:
$$ 
\check{H}^1(M,EH) \to \check{H}^1(M, H\to G\ltimes EH) 
\to \check{H}^1(M, H\to G) 
$$ 
Since $EH$ is contractible and $M$ is paracompact Hausdorff,
$\check{H}^1(M,EH)$ is easily seen to be trivial, so the 
map $\check{H}^1(M,H\to G\ltimes EH)\to \check{H}^1(M,H\to G)$ 
is injective.  To see that this map is surjective, note that
there is a homomorphism of crossed modules going back:
\[
\xymatrix{ 
1 \ar[r] & 1 \ar[r] \ar[d] & H \ar[d] \ar[r]^-1 & H \ar[d]^-{t} 
\ar@/_1pc/[l]_-1 \ar[r] & 1 \\ 
1 \ar[r] & EH \ar[r] & G \ltimes EH \ar[r] & G \ar[r] 
\ar@/_1pc/[l]_-i & 1 } 
\]
where $i$ is the natural inclusion of $G$ in the
semidirect product $G\ltimes EH$.  This homomorphism going
back `splits' our exact sequence of crossed modules.
It follows that $\check{H}^1(M,H\to G \ltimes EH)\to
\check{H}^1(M, H\to G)$ is onto, so we have a bijection
$$ 
\check{H}^1(M, H\to G\ltimes EH) \cong \check{H}^1(M,H\to G) =
\check{H}^1(M, \G)
$$
completing the proof.

\subsection{Remarks on Theorem \ref{classifying}}
\label{remarks}

Theorem \ref{classifying}, asserting the existence of a classifying space 
for first \v{C}ech cohomology with coefficients in a topological 2-group, 
was originally stated in a preprint by Jur\v{c}o \cite{Jurco}.  However,
the argument given there was missing some details.  In essence, the
Jur\v{c}o's argument boils down to the following: he constructs
a map $\check{H}^1(M,|\G|)\to \check{H}^1(M,\G)$ and sketches the
construction of a map $\check{H}^1(M,\G)\to \check{H}^1(M,|\G|)$.  The
construction of the latter map however requires some further
justification: for instance, it is not obvious that one can choose a
classifying map satisfying the cocycle property listed on the top of
page 13 of \cite{Jurco}.  Apart from this, it is not demonstrated that
these two maps are inverses of each other.

As mentioned earlier, Jur\v{c}o and also Baas,
B\"okstedt and Kro \cite{BBK} use a different approach to construct
a classifying space for a topological 2-group $\G$.  In their
approach, $\G$ is regarded as a topological 2-groupoid with one
object.  There is a well-known nerve construction that turns any
2-groupoid (or even any 2-category) into a simplicial set \cite{Duskin}.
Internalizing this construction, these authors turn the topological
2-groupoid $\G$ into a simplicial space, and then take the geometric
realization of that to obtain a space.  Let us denote this space by $B\G$.  
This is the classifying space used by Jur\v{c}o and Baas--B\"okstedt--Kro.  
It should be noted that the assumption that $\G$ is a well-pointed  
2-group ensures that the nerve of the 2-groupoid $\G$ is a `good' 
simplicial space in the sense of Segal; this `goodness' condition is 
important in the work of Baas, B\"okstedt and Kro \cite{BBK}.

Baas, B\"okstedt and Kro also consider a third way to construct
a classifying space for $\G$.  If we take the nerve $N\G$ of $\G$ we get a 
simplicial group, as described in Section~\ref{2groups} above.  
By thinking of each group of $p$-simplices $(N\G)_p$ as a groupoid with 
one object, we can think of $N\G$ as a simplicial groupoid.  
From $N\G$ we can obtain a bisimplicial space $NN\G$ by applying the 
nerve construction to each groupoid $(N\G)_p$.  $NN\G$ is sometimes
called the `double nerve', since we apply the nerve construction twice.  
From this bisimplicial space $NN\G$ we can form an ordinary simplicial 
space $dNN\G$ by taking the diagonal.  Taking the geometric realization 
of this simplicial space, we obtain a space $|dNN\G|$.

It turns out that this space $|dNN\G|$ is homeomorphic to 
$B|\G|$ \cite{Braho, Quillen}.  It can also be shown that the spaces 
$|dNN\G|$ and $B\G$ are homotopy equivalent --- but although this fact 
seems well-known to experts, we have been unable to find a reference 
in the case of a \emph{topological} 2-group $\G$.  For ordinary 2-groups
(without topology) the relation between all three nerves was worked out by 
Moerdijk and Svensson \cite{MS} and Bullejos and Cegarra \cite{BC}.  
In any case, since we do not use these facts in our arguments, we 
forgo providing the proofs here.

\subsection{Proof of Lemma \ref{lem1}} 
\label{proof.lem.1}

Suppose $\G$ is a well-pointed topological 2-group with 
topological crossed module $(G,H,t,\alpha)$, and let 
$|\G|$ be the geometric realization of its nerve.
We shall prove that there is a topological 2-group $\hat{\G}$ 
fitting into a short exact sequence of topological 2-groups 
\begin{equation} 
\label{eq: exact seq of 2-groups} 
1\to H\to \hat{\G}\to \G\to 1
\end{equation}
where $H$ is the discrete topological $2$-group associated to the 
topological group $H$.  On taking nerves and then geometric 
realizations, this gives an exact sequence of groups: 
$$ 
1\to H\to |\hat{\G}| \to |\G|\to 1
$$ 
Redescribing the 2-group $\hat{\G}$ with the help of some work by Segal,
we shall show that $|\hat{\G}| \cong G\ltimes EH$ and thus $|\G| \cong
(G \ltimes EH)/H$.  Then we prove that the above sequence is an exact
sequence of {\em topological} groups: this requires checking that
$|\hat{\G}| \to |\G|$ is a Hurewicz fibration.  We conclude by showing that
$|\hat{\G}|$ is well-pointed.

To build the exact sequence of $2$-groups in equation 
\eqref{eq: exact seq of 2-groups}, we construct the corresponding
exact sequence of topological crossed modules.  This
takes the following form:
$$ 
\xymatrix{ 
1 \ar[r] & 1 \ar[d] \ar[r] & H \ar[d]^-{t'} \ar[r]^-1 & H \ar[d]^-{t} \ar[r] & 1  \\ 
1\ar[r] & H \ar[r]^-{f} & G\ltimes H \ar[r]^-{f'} & G \ar[r] & 1} 
$$ 
Here the crossed module $(G\ltimes H,H,t',\alpha')$ is 
defined as follows: 
\begin{align*} 
& t'(h) = (1,h) \\ 
& \a'(g,h)(h') = \a(t(h)g)(h') 
\end{align*} 
while $f$ and $f'$ are given by
\[
\begin{array}{rcl}
f\maps  H &\to& G\ltimes H \\ 
 h &\mapsto & (t(h),h^{-1}) \\ 
\\ 
f'\maps G\ltimes H &\to& G \\ 
 (g,h) &\mapsto& t(h)g 
\end{array} 
\]
It is easy to check that these formulas define an exact sequence of
topological crossed modules.  The corresponding exact sequence 
of topological $2$-groups is
$$ 
1\to H\to \hat{\G} \to \G\to 1 
$$ 
where $\hat{\G}$ 
denotes the topological $2$-group associated to the
topological crossed module $(G\ltimes H,H,t',\alpha')$.  

In more detail, the $2$-group $\hat{\G}$ has 
\begin{align*} 
& \Ob(\hat{\G}) = G\ltimes H \\ 
& \Mor(\hat{\G}) = (G\ltimes H)\ltimes H \\ 
& s((g,h),h') = (g,h), \qquad t((g,h),h') = (g,h'h) \\ 
& i(g,h) = ((g,h),1), \qquad ((g,h'h),h'') \circ ((g,h),h') = ((g,h),h''h')
\end{align*} 
Note that there is an isomorphism 
$(G\ltimes H)\ltimes H \cong G\ltimes H^2$ sending
$((g,h),h')$ to $(g,(h,h'h))$.  Here by $G\ltimes H^2$ we mean the 
semidirect product 
formed with the diagonal action of $G$ on $H^2$, namely
$g(h,h') = (\a(g)(h),\a(g)(h'))$.  Thus the group 
$\Mor(\hat{\G})$ is isomorphic to $G\ltimes H^2$.  

We can give a clearer description of the 2-group $\hat{\G}$ using the work 
of Segal \cite{Segal1}.  Segal noted that for any topological group $H$,
there is a $2$-group $\overline{H}$ with one object for each element of 
$H$, and one morphism from any object to any other.  In other words, 
$\overline{H}$ is the 2-group with:
\begin{align*} 
& \Ob(\overline{H}) = H \\ 
& \Mor(\overline{H}) = H^2 \\ 
& s(h,h') = h, \qquad t(h,h') = h' \\ 
& i(h) = (h,h), \qquad (h',h'') \circ (h,h') = (h,h'')
\end{align*} 
Moreover, Segal proved that the geometric realization $|\overline{H}|$ 
of the nerve of $\overline{H}$ is a model for $EH$.  Since 
$G$ acts on $H$ by automorphisms, we can define a `semidirect product' 
2-group $G\ltimes \overline{H}$ with 
\begin{align*} 
& \Ob(G\ltimes \overline{H}) = G\ltimes H \\ 
& \Mor(G\ltimes \overline{H}) = G\ltimes H^2 \\ 
& s(g,(h,h')) = (g,h),\qquad t(g,(h,h')) = (g,h') \\ 
& i(g,h) = (g,(h,h)),\qquad (g,(h',h'')) \circ (g,(h,h')) = (g,(h,h''))
\end{align*} 
The isomorphism $(G\ltimes H)\ltimes H \cong G\ltimes H^2$ above can 
then be interpreted as an isomorphism $\Mor(\hat{\G})\cong \Mor(G\ltimes 
\overline{H})$.  It is easy to check that this isomorphism  
is compatible with the structure maps for $\hat{\G}$ and $G\ltimes 
\overline{H}$, so we have an isomorphism of topological 2-groups:
$$ 
\hat{\G} \cong G\ltimes \overline{H}
$$ 
It follows that the nerve $N\hat{\G}$ of $\hat{\G}$ is isomorphic as a 
simplicial topological group to the nerve of $G\ltimes \overline{H}$.  
As a simplicial space it is clear that $N(G\ltimes \overline{H}) = 
G\times N\overline{H}$. We need to identify the simplicial group structure 
on $G\times N\overline{H}$.  

From the definition of the products on $\Ob(G\ltimes \overline{H})$ and 
$\Mor(G\ltimes \overline{H})$, it is clear that the product on $N(G\ltimes 
\overline{H})$ is given by the simplicial map 
\[   (G\times N\overline{H})\times (G\times N\overline{H})\to 
G\times N\overline{H} \]
defined on $p$-simplices by 
$$ 
\left((g,(h_1,\ldots, h_p)),(g',(h_1',\ldots, h'_p))\right) 
\mapsto (gg',(h_1\a(g)(h'_1),\ldots, h_p\a(g)(h'_p))) 
$$
Thus one might well call $N(G\ltimes \overline{H})$ the `semidirect 
product' $G\ltimes N\overline{H}$.  Since geometric realization preserves 
products, it follows that there is an isomorphism of 
topological groups 
$$ 
|\hat{\G}| \cong G\ltimes EH.  
$$ 
Here the semidirect product is formed using the action of $G$ on $EH$ 
induced from the action of $G$ on $H$. 
Finally note that $H$ is embedded as a normal subgroup of $G\ltimes EH$ 
through 
\[ 
\begin{array}{ccc}
H &\to& G\ltimes EH  \\
h &\mapsto& (t(h),h^{-1}) 
\end{array}
\]
It follows that the exact sequence of groups 
$1\to H\to |\hat{\G}|\to |\G|\to 1$ 
can be identified with 
\begin{equation} 
\label{eq: fundamental extension} 
1\to H\to G\ltimes EH\to |\G|\to 1
\end{equation}
It follows that $|\G|$ is isomorphic to the 
quotient $G\ltimes_H EH$ of $G\ltimes EH$ by the normal subgroup $H$.  
This amounts to factoring out by the action of $H$ on $G\ltimes EH$ 
given by $h(g,x) = (t(h)g,xh^{-1})$.  

Next we need to show that equation \eqref{eq: fundamental extension} 
specifies an exact sequence of {\em topological} groups: in particular, 
that the map $G\ltimes EH \to |\G| = G\ltimes_H EH$ is a Hurewicz fibration.   
To do this, we prove that the following diagram is a pullback: 
$$ 
\xymatrix{ 
G\ltimes EH\ar[d] \ar[r] & EH \ar[d] \\ 
G\ltimes_H EH \ar[r] & BH } 
$$ 
Since $H$ is well pointed, $EH\to BH$ is a numerable principal
bundle (and hence a Hurewicz fibration) by the results of May
\cite{May1} referred to earlier.  The statement above now follows, as
Hurewicz fibrations are preserved under pullbacks.

To show the above diagram is a pullback, we construct a homeomorphism  
$$ 
\a\colon (G\ltimes_H EH) \times_{BH} EH \to G\ltimes EH 
$$ 
whose inverse is the canonical map $\b\colon G\ltimes EH\to (G\ltimes_H EH) 
\times_{BH} EH$.  To do this, suppose that 
$([g,x],y)\in (G\ltimes_H EH) \times_{BH} EH$.  Then 
$x$ and $y$ belong to the same fiber of $EH$ over $BH$, so 
$y^{-1}x\in H$.  We set 
$$ 
\a([g,x],y) = ( t(y^{-1}x)g,y)
$$ 
A straightforward calculation shows that $\a$ is well defined
and that $\a$ and $\b$ are inverse to one another.

To conclude, we need to show that $|\G|$ is a well-pointed 
topological group.  For this it is sufficient to show that $N\G$ is a 
`proper' simplicial space in the sense of May \cite{May2}
(note that we can replace his `strong' NDR pairs with NDR pairs).
For, if we follow May and denote by $F_p|\G|$ the image of 
$\coprod^p_{i=0} \Delta^i\times N\G_i$ in $|\G|$, it then follows from 
his Lemma 11.3 that $(|\G|,F_p|\G|)$ is an NDR pair for all $p$.  
In particular $(|\G|,F_0|\G|)$ is an NDR pair.  Since $F_0|\G| = G$ 
and $(G,1)$ is an NDR pair, it follows that $(|\G|,1)$ is an NDR pair:
that is, $|\G|$ is well pointed.    

We still need to show that $N\G$ is proper.  In fact it suffices to
show that $N\G$ is a `good' simplicial space in the sense of Segal
\cite{Segal3}, meaning that all the degeneracies $s_i\colon
N\G_n\to N\G_{n+1}$ are closed cofibrations.  The reason for this is
that every good simplicial space is automatically proper --- see the
proof of Lewis' Corollary 2.4(b) \cite{GaunceLewis}.  To see that
$N\G$ is good, note that every degeneracy homomorphism $s_i\colon
N\G_n\to N\G_{n+1}$ is a section of the corresponding face
homomorphism $d_i$, so $N\G_{n+1}$ splits as a semidirect
product $N\G_{n+1}\cong N\G_n\ltimes \ker(d_i)$.  Therefore, $s_i$ is
a closed cofibration provided that $\ker(d_i)$ is well pointed.  But
$\ker(d_i)$ is a retract of $N\G_{n+1}$, so $\ker(d_i)$ will
be well pointed if $N\G_{n+1}$ is well pointed.  For this, note that
$N\G_{n+1}$ is isomorphic as a space to $G\times
H^{n+1}$.  Since the groups $G$ and $H$ are well pointed by
hypothesis, it follows that $N\G_{n+1}$ is well pointed.  Here we have
used the fact that if $X\to Y$ and $X'\to Y'$ are closed cofibrations
then $X\times X'\to Y\times Y'$ is a closed cofibration.  \qed
  
\subsection{Proof of Lemma \ref{lem2}} 
\label{proof.lem.2}

Suppose that $M$ is a topological space admitting good covers.
Also suppose that $1\to H\stackrel{t}{\to} G
\stackrel{p}{\to} K\to 1$ is an exact sequence of topological groups.

This data gives rise to a topological crossed module 
$H \stackrel{t}{\to} G$ where $G$ acts on $H$ by conjugation.  
For short we denote this by $H \to G$.  
The same data also gives a topological crossed module $1 \to K$.  There is 
a homomorphism of
crossed modules from $H \to G$ to $1 \to K$, arising from this commuting 
square:
\[ 
\xymatrix{ 
H \ar[r] \ar[d]^{t} & 1 \ar[d] \\ 
G \ar[r]^{p} & K } 
\]
Call this homomorphism $\alpha$.  It yields a map
\[ \alpha_\ast \maps \check{H}^1(M,H \to G) \to \check{H}^1(M,1 \to K). \]
Note that $\check{H}^1(M,1 \to K)$ is just the ordinary \v{C}ech 
cohomology $\check{H}^1(M,K)$.  
To prove Lemma \ref{lem2}, we need to construct an inverse
\[ \beta \maps \check{H}^1(M,K) \to \check{H}^1(M,H \to G). \]
Let $\cU = \{U_i\}$ be a good cover of $M$; then, as noted in
Section~\ref{nonabelian} there is a bijection
\[ \check{H}^1(M,K) = \check{H}^1(\cU,K) \]
Hence to define the map $\beta$ it is sufficient to define a map
$\beta\colon \check{H}^1(\cU,K)\to \check{H}^1(\cU,H\to G)$.  Let
$k_{ij}$ be a $K$-valued \v{C}ech 1-cocycle subordinate to $\cU$.
Then from it we construct a \v{C}ech 1-cocycle $(g_{ij},h_{ijk})$
taking values in $H \to G$ as follows.  Since the spaces $U_i\cap U_j$
are contractible and $p \maps G \to K$ is a Hurewicz fibration, 
we can lift the maps $k_{ij} \maps U_i\cap U_j \to K$ to maps
$g_{ij}\maps U_i\cap U_j\to G$.  The $g_{ij}$ need not satisfy the
cocycle condition for ordinary \v{C}ech cohomology, but instead we
have
\[
t(h_{ijk}) g_{ij}g_{jk} = g_{ik}
\]
for some unique $h_{ijk}\colon U_i\cap U_j\cap U_k\to H$.   In terms
of diagrams, this means we have triangles
\[
\xy 
(0,0)*+{\bullet}="1"; 
(-12,-20.78)*+{\bullet}="2"; 
(12,-20.78)*+{\bullet}="3"; 
{\ar^-{g_{ij}} "2";"1"}; 
{\ar^-{g_{jk}} "1";"3"}; 
{\ar_-{g_{ik}} "2";"3"}; 
{\ar@2{->}_-{h_{ijk}} (0,-7)*{};(0,-18)*{}};
\endxy
\]
The uniqueness of $h_{ijk}$ follows from the fact that the 
homomorphism $t\colon H\to G$ is injective.
To show that the pair $(g_{ij},h_{ijk})$ defines a \v{C}ech cocycle 
we need to check that the tetrahedron~\eqref{tetrahedron} commutes.  
However, this follows from the commutativity of the corresponding
tetrahedron built from triangles of this form:
\[
\xy 
(0,0)*+{\bullet}="1"; 
(-12,-20.78)*+{\bullet}="2"; 
(12,-20.78)*+{\bullet}="3"; 
{\ar^-{k_{ij}} "2";"1"}; 
{\ar^-{k_{jk}} "1";"3"}; 
{\ar_-{k_{ik}} "2";"3"}; 
{\ar@2{->}_-{1} (0,-7)*{};(0,-18)*{}};
\endxy
\]
and the injectivity of $t$.  

Let us show that this construction gives a well-defined map
\[  \beta \maps \check{H}^1(M,K) = 
\check{H}^1(\cU,K) \to \check{H}^1(M,H\to G) \] 
sending $[k_{ij}]$ to $[g_{ij}, h_{ijk}]$.  
Suppose that $k_{ij}'$ is another $K$-valued \v{C}ech 1-cocycle subordinate to 
$\cU$, such that $k'_{ij}$ and $k_{ij}$ are cohomologous.  Starting from the 
cocycle $k_{ij}'$ we can construct (in the same manner as above) a cocycle 
$(g'_{ij},h'_{ijk})$ taking values in $H\to G$.  Our task is to show that 
$(g_{ij},h_{ijk})$ and $(g'_{ij},h'_{ijk})$ are cohomologous.  Since 
$k_{ij}$ and $k'_{ij}$ are cohomologous there exists a 
family of maps 
$\kappa_i\colon U_i\to K$ 
fitting into the naturality square 
\[ 
\xymatrix{ 
\bullet \ar[r]^-{k_{ij}(x)} \ar[d]_-{\kappa_i(x)} & \bullet \ar[d]^-{\kappa_j(x)} \\ 
\bullet \ar[r]_-{k_{ij}'(x)} & \bullet }
\] 
Choose lifts $f_i\colon U_i\to G$ of the various $\kappa_i$.  Since $p(g_{ij}) = p(g_{ij}') = k_{ij}$ 
and $p(f_i) = \kappa_i$, $p(f_j) = \kappa_j$ there is a unique map $\eta_{ij}\maps U_{i}\cap U_{j} \to H$ 
\[ 
t(\eta_{ij}) g_{ij}f_j = f_ig_{ij}' .\]
So, in terms of diagrams, we have the following squares:

\[ 
\xy 
(-7.5,0)*+{\bullet} ="1"; 
(7.5,0)*+{\bullet} = "2"; 
(-7.5,-15)*+{\bullet} = "3"; 
(7.5,-15)*+{\bullet}="4"; 
{\ar^-{g_{ij}} "1";"2"}; 
{\ar_-{f_i} "1";"3"}; 
{\ar^-{f_j} "2";"4"}; 
{\ar_-{g'_{ij}} "3";"4"}; 
{\ar@2{->}_-{\eta_{ij}} (6,-2)*{};(-4,-13)*{}};
\endxy 
\]
The triangles and squares defined so far fit together to form prisms:
\begin{center}
\begin{picture}(200,270)
\includegraphics{prism.eps}
  \put(-1,100){$f_k$}
  \put(-170,100){$f_i$}
  \put(-80,160){${}_{f_j}$}
  \put(-90,-7){$g'_{ik}$}
  \put(-100,187){$g_{ik}$}
  \put(-128,43){${}_{g'_{ij}}$}
  \put(-128,237){${}_{g_{ij}}$}
  \put(-47,43){${}_{g'_{jk}}$}
  \put(-47,237){${}_{g_{jk}}$}
  \put(-95,16){${}_{h'_{ijk}}$}
  \put(-95,212){${}_{h_{ijk}}$}
  \put(-95,115){$\eta_{ik}$}
  \put(-132,142){${}_{\eta_{ij}}$}
  \put(-35,121){${}_{\eta_{jk}}$}
\end{picture}
\end{center}
It follows from the injectivity of the homomorphism $t$ that 
these prisms commute, and therefore that $(g_{ij},h_{ijk})$
and $(g'_{ij},h'_{ijk})$ are cohomologous. Therefore we have a 
well-defined map $\check{H}^1(\cU,K)\to \check{H}^1(\cU,H\to G)$ 
and hence a well-defined map $\beta\colon \check{H}^1(M,K)\to 
\check{H}^1(M,H\to G)$.    

Finally we need to check that $\alpha$ and $\beta$ are inverse to one
another.  It is obvious that $\alpha\circ \beta$ is the identity on
$\check{H}^1(M,K)$.  To see that $\beta \circ \alpha$ is the identity
on $\check{H}^1(M,H\to G)$ we argue as follows.  Choose a cocycle
$(g_{ij},h_{ijk})$ subordinate to a good cover $\cU = \{U_i\}$.  
Then under $\alpha$ the cocycle $[g_{ij},h_{ijk}]$ is sent to the 
$K$-valued cocycle $[p(g_{ij})]$.  But then we may take $g_{ij}$ as 
our lift of $p(g_{ij})$ in the definition of $\beta(p(g_{ij}))$.  
It is then clear that $(\beta \circ \alpha) [g_{ij},h_{ijk}] = 
[g_{ij},h_{ijk}]$.  \qed

\vskip 1em
At this point a remark is in order.  The proof of the above lemma is 
one place where the definition of $\check{H}^1(M,H\to G)$ in terms 
of hypercovers would lead to simplifications, and would allow us to 
replace the hypothesis that the map underlying the homomorphism $G\to K$ 
was a fibration with a less restrictive condition.  The homomorphism 
of crossed modules 
\[ 
\xymatrix{ 
H \ar[r] \ar[d]^{t} & 1 \ar[d] \\ 
G \ar[r]^{p} & K } 
\]
gives a homomorphism between the associated 2-groups and
hence a simplicial map between the nerves of the associated
2-groupoids.  It turns out that this simplicial map belongs to a
certain class of simplicial maps with respect to which a subcategory
of simplicial spaces is localized.  In the formalism of hypercovers, 
for $M$ paracompact, 
the nonabelian cohomology $\check{H}^1(M,\G)$ with coefficients in a 2-group
$\G$ is defined as a certain set of morphisms in this localized
subcategory.  It is then easy to see that the induced map
$\check{H}^1(M,H\to G)\to \check{H}^1(M,K)$ is a bijection.

\subsection{Proof of Lemma \ref{lem3}} 
\label{proof.lem.3}

Suppose that 
$$ 1 \to \G_0 \stackrel{f}{\to} \G_1 \stackrel{p}{\to} \G_2 \to 1$$ 
is a short exact sequence of topological 2-groups, so that we have a short
exact sequence of topological crossed modules:
$$ 
\xymatrix{ 
1 \ar[d] \ar[r] & H_0 \ar[d]_-{t_0} \ar[r]^-{f} & H_1 \ar[d]_-{t_1} \ar[r]^-{p}  
& H_2 \ar[d]_-{t_2} \ar[r] & 1 \ar[d] \\ 
1\ar[r] & G_0 \ar[r]^-f & G_1 \ar[r]^-{p} & G_2 \ar[r] & 1} 
$$ 
Also suppose that $\cU = \{U_i\}$ is a good cover of $M$, and that
$(g_{ij},h_{ijk})$ is a cocycle representing a class in 
$\check{H}^1(\cU,\G_1)$.  We claim that
the image of 
\[ f_\ast \maps \check{H}^1(M,\G_0) \to \check{H}^1(M,\G_1)  \]
equals the kernel of 
\[ p_\ast \maps \check{H}^1(M,\G_1) \to \check{H}^1(M,\G_2) . \]
If the class $[g_{ij},h_{ijk}]$ is in the image of $f_\ast$,
it is clearly in the kernel of $p_\ast$.  Conversely, 
suppose it is in kernel of $p_\ast$.  We need to show that it
is in the image of $f_\ast$.

The pair $(p(g_{ij}),p(h_{ijk}))$ is cohomologous to the trivial
cocycle, at least after refining the cover $\cU$, so there exist
$x_i\maps U_i\to G_2$ and $\xi_{ij} \maps U_i\cap U_j \to H_2$ such
that this diagram commutes:
\begin{center}
\begin{picture}(180,270)
\includegraphics{prism.eps}
  \put(-1,100){$x_k$}
  \put(-170,100){$x_i$}
  \put(-80,160){${}_{x_j}$}
  \put(-85,-7){$1$}
  \put(-100,187){$p(g_{ik})$}
  \put(-128,43){${}_1$}
  \put(-129,237){${}_{p(g_{ij})}$}
  \put(-47,43){${}_1$}
  \put(-47,237){${}_{p(g_{jk})}$}
  \put(-95,16){${}_1$}
  \put(-108,212){${}_{p(h_{ijk})}$}
  \put(-95,115){$\xi_{ik}$}
  \put(-132,142){${}_{\xi_{ij}}$}
  \put(-35,121){${}_{\xi_{jk}}$}
\end{picture}
\end{center}

Since $p \maps G_1 \to G_2$ is a fibration and 
$U_i$ is contractible, we can lift $x_i$ to
a map $\hat{x}_i \maps U_i \to G_1$.  Similarly,
we can lift $\xi_{ij}$ to a map $\hat{\xi}_{ij} \maps U_i\cap U_j \to H_1$.   
There are then unique maps $\gamma_{ij} \maps U_i\cap U_j \to G_1$ giving 
squares like this:
\[
\xy
(-7.5,0)*+{\bullet} ="1";
(7.5,0)*+{\bullet} = "2";
(-7.5,-15)*+{\bullet} = "3";
(7.5,-15)*+{\bullet}="4";
{\ar^-{g_{ij}} "1";"2"};
{\ar_-{\hat{x}_i} "1";"3"};
{\ar^-{\hat{x}_j} "2";"4"};
{\ar_-{\c_{ij}} "3";"4"};
{\ar@2{->}_-{\hat{\xi}_{ij}} (6,-2)*{};(-4,-13)*{}};
\endxy
\]
namely 
\[ 
 \c_{ij} = \hat{x}_ig_{ij}\hat{x}_j^{-1} t(\hat{\xi}_{ij}) 
\]
Similarly, there are unique maps $c_{ijk} \maps U_i\cap U_j\cap U_k \to H_1$
making this prism commute:
\begin{center}
\begin{picture}(180,270)
\includegraphics{prism.eps}
  \put(-1,100){$\hat{x}_k$}
  \put(-170,100){$\hat{x}_i$}
  \put(-80,160){${}_{\hat{x}_j}$}
  \put(-90,-7){$\c_{ik}$}
  \put(-100,187){$g_{ik}$}
  \put(-128,43){${}_{\c_{ij}}$}
  \put(-129,237){${}_{g_{ij}}$}
  \put(-47,43){${}_{\c_{jk}}$}
  \put(-47,237){${}_{g_{jk}}$}
  \put(-95,16){${}_{c_{ijk}}$}
  \put(-95,212){${}_{h_{ijk}}$}
  \put(-95,115){$\hat{\xi}_{ik}$}
  \put(-132,142){${}_{\hat{\xi}_{ij}}$}
  \put(-35,121){${}_{\hat{\xi}_{jk}}$}
\end{picture}
\end{center}
To define $c_{ijk}$, we simply compose the 2-morphisms on the
sides and top of the prism.

Applying $p$ to the prism above we obtain the previous prism.
So, $\c_{ij}$ and $c_{ijk}$ must take values in the kernel of 
$p\maps G_1 \to G_2$ and $p \maps H_1 \to H_2$, respectively.
It follows that $\c_{ij}$ and $c_{ijk}$ take values in the image
of $f$. 

The above prism says that $(\c_{ij},c_{ijk})$ is cohomologous
to $(g_{ij}, h_{ijk})$, and therefore a cocycle in its own
right.  Since $\c_{ij}$ and $c_{ijk}$ take values in the image
of $f$, they represent a class in the image of 
\[ f_\ast \maps \check{H}^1(M,\G_0) \to \check{H}^1(M,\G_1) . \]
So, $[g_{ij}, h_{ijk}] = [\c_{ij},c_{ijk}]$ is in the image of
$f_\ast$, as was to be shown.   \qed

\subsection*{Proof of Theorem \ref{string}}
\label{proof.thm.2}

The following proof was first described to us by Matt Ando \cite{Ando}, and
later discussed by Greg Ginot \cite{Ginot}.

Suppose that $G$ is a simply-connected, compact, simple Lie group.  
Then the string group $\hat{G}$ of $G$ fits into a short exact sequence 
of topological groups 
$$ 
1\to K(\ZZ,2)\to \hat{G}\to G\to 1 
$$ 
for some realization of the Eilenberg-Mac Lane space $K(\ZZ,2)$ 
as a topological group.  Applying the classifying space functor $B$ 
to this short exact sequence gives rise to a fibration 
$$ 
K(\ZZ,3)\to B\hat{G}\stackrel{p}{\to} BG. 
$$ 
We want to compute the rational cohomology of $B\hat{G}$.  




We can use the Serre spectral sequence to compute $H^*(B\hat{G},\QQ)$.  
Since $BG$ is simply connected the $E_2$ term of this spectral sequence is 
$$ 
E_2^{p,q} = H^p(BG,\QQ)\otimes H^q(K(\ZZ,3),\QQ). 
$$ 
Because $K(\ZZ,3)$ is rationally indistinguishable from $S^3$, the first 
nonzero differential is $d_4$.  Furthermore, the differentials of this 
spectral sequence are all derivations.  It follows that $d_4(y\otimes x_3) = 
(-1)^{p}y\otimes d_4(x_3)$ if $y\in H^p(BG,\QQ)$.  It is not hard to identify 
$d_4(x_3)$ with $c$, the class in $H^4(BG,\QQ)$ which is the transgression 
of the generator $\nu$ of $H^3(G,\QQ) = \QQ$.  It follows that the spectral 
sequence collapses at the $E_5$ stage with 
$$ 
E_5^{p,q} = E_{\infty}^{p,q} = \begin{cases} 
0 & \text{if}\ q > 0 \\ 
H^p(BG,\QQ)/\langle c\rangle & \text{if}\ q = 0. 
\end{cases} 
$$ 
One checks that all the subcomplexes $F^iH^*(B\hat{G},\QQ)$ in the 
filtration of $H^*(B\hat{G},\QQ)$ are zero for $i\geq 1$.  Hence 
$H^p(B\hat{G},\QQ) = E^{p,0}_{\infty} = H^p(BG,\QQ)/\langle c \rangle$ 
and so Theorem \ref{string} is proved.

\subsubsection*{Acknowledgements}  We would like to thank Matthew
Ando and Greg Ginot for showing us the proof of Theorem \ref{string}, 
Peter May for helpful remarks on classifying spaces, and Toby Bartels,
Branislav Jur\v{c}o, and Urs Schreiber for many discussions on higher 
gauge theory.  We also thank Nils Baas, Bj\o rn Dundas, Eric Friedlander, 
Bj\o rn Jahren, John Rognes, Stefan Schwede, and Graeme Segal for
organizing the 2007 Abel Symposium and for useful discussions.
Finally, DS has received support from Collaborative Research 
Center 676: `Particles, Strings, and the Early Universe'.

\end{document}